\documentclass[11pt]{amsart}
\usepackage{amssymb, amsmath, amsthm}
\usepackage[latin1]{inputenc}
\usepackage{graphicx}
\usepackage[final]{hyperref}
\usepackage{color}
\usepackage{algorithm}
\usepackage{algorithmicx}
\usepackage{algpseudocode}

\usepackage[a4paper, centering]{geometry}
\usepackage{color}
\usepackage{graphicx}
\graphicspath{{./pics/}}
\newtheorem{theorem}{Theorem}
\newtheorem{proposition}[theorem]{Proposition}

\newtheorem{corollary}[theorem]{Corollary}
\theoremstyle{definition}

\theoremstyle{remark}

\newtheorem{conjecture}[theorem]{Conjecture}

\def\R{\mathbb{R}}

\def\N{\mathbb{N}}

\definecolor{verde}{RGB}{20,150,100}

\newcommand{\Om}{\Omega}

\newcommand{\ra}{\rightarrow}

\def \d{\delta}

\newcommand{\vps}{\varepsilon}

\newcommand{\HN}{{\mathcal H}^{N-1}}
\newcommand{\sm}{\setminus}
\newcommand{\lra}{\longrightarrow}
\newcommand{\sr}{\stackrel}
\newcommand{\sq}{\subseteq}

\def\G{\Gamma}

\begin{document}

\title[]{Phase field approach to optimal packing problems \\ and related Cheeger clusters }
\thanks{This work  was supported by the  ANR-15-CE40-0006 COMEDIC project and GNAMPA (INDAM)}

\author[]{Beniamin Bogosel, Dorin Bucur, Ilaria Fragal\`a}
\thanks{}

\address [Beniamin Bogosel]{CMAP UMR 7641 École Polytechnique CNRS, Route de Saclay, 91128 Palaiseau Cedex (France)}
\email {beniamin.bogosel@cmap.polytechnique.fr}

\address[Dorin Bucur]{
Laboratoire de Math\'ematiques UMR 5127 \\
Universit\'e de Savoie,  Campus Scientifique \\
73376 Le-Bourget-Du-Lac (France)
}
\email{dorin.bucur@univ-savoie.fr}

\address[Ilaria Fragal\`a]{
Dipartimento di Matematica \\ Politecnico  di Milano \\
Piazza Leonardo da Vinci, 32 \\
20133 Milano (Italy)
}
\email{ilaria.fragala@polimi.it}


\date{\today}
\maketitle

\begin{abstract}  In a fixed domain of $\R ^N$  we study the asymptotic behaviour of optimal clusters  associated to $\alpha$-Cheeger constants and natural energies like the sum or maximum: we prove that, as the parameter $\alpha$ converges to the ``critical" value $\Big (\frac{N-1}{N}\Big ) _+$, optimal Cheeger clusters converge to solutions of  different packing problems for balls, depending on the energy under consideration. As well, we propose an efficient phase field approach based on a multiphase Gamma convergence result of Modica-Mortola type, in order to compute $\alpha$-Cheeger constants, optimal clusters and, as a consequence of the asymptotic result, optimal packings. Numerical experiments are carried over in two and three space dimensions.
  \end{abstract}



\section{Introduction and statement of the results}

Let $N \ge 2$ be the space dimension and $\alpha > \frac{N-1}{N}$ be a fixed constant. For every bounded measurable subset  $E$ of $\R^N$, we define its $\alpha$-Cheeger constant by
\begin{equation}\label{halpha}
h_\alpha (E):= \min \{\frac{\HN (\partial ^* \Omega)}{|\Om|^\alpha}\  : \ \Om \sq E, \ \Om \mbox{ measurable}\}.
\end{equation}
Above, $\HN$ denotes the $(N-1)$-dimensional Hausdorff measure and, if $\Om^*$ has finite perimeter,  $\partial ^* \Omega$ is its reduced boundary, in the measure theoretical sense. If $\Om$ is the empty set or it has positive measure but does not have finite perimeter, the ratio $\frac{\HN (\partial ^* \Omega)}{|\Om|^\alpha}$ is assumed by convention to be equal to $+\infty$.

For $\alpha=1$, definition \eqref{halpha} corresponds to the classical Cheeger constant, which  was thoroughly studied in the last years, 
see for instance the review papers \cite{Leo, Pa};  for $\alpha\neq 1$, and strictly larger than the scale invariance exponent $\frac{N-1}{N}$, the notion of $\alpha$-Cheeger constant is a variant which has appeared in the literature more recently, we refer in particular to \cite{PS17} and references therein.

Object of this paper are the following optimal partition problems

\begin{equation}\label{f:pb1.b}\min  \Big \{\max_{i=1, \dots, k }  h_\alpha(E_i)  :\   (E_1, \dots, E_k) \in \mathcal P _k (D)  \Big \}
\end{equation}

\begin{equation}\label{f:pb2.b}\min  \Big \{\sum_{i=1} ^ k    h_\alpha(E_i)   :\  (E_1, \dots, E_k) \in \mathcal P _k (D)  \Big \}\,,
\end{equation}
where $D$ is a given open bounded Lipschitz subset of  $\R^N$, and 
$$\mathcal P _k (D)=\Big \{ (E_1, \dots, E_k) \ : \ \forall i,j =1, \dots, k, \ E_i \sq D, \ E_i \cap E_j =\emptyset, \ E_i \mbox{ measurable} \Big \}.$$

Notice that the condition that the union of the sets $E_i$ covers  the given box $D$ (up to a negligible set) is not required in the definition $\mathcal P _k (D)$, 
so that there is some abuse of notation in adopting the usual  epithet of optimal ``partitions'' as done above, and in the sequel we prefer to speak rather about ``clusters''. 
Thus, solutions of \eqref{f:pb1.b}-\eqref{f:pb2.b} will be generically called {\it $\alpha$-Cheeger clusters}.

For $\alpha=1$, problem \eqref{f:pb2.b} has been firstly studied by Caroccia in the paper \cite{Car17}, where the existence of solutions and some regularity results for the free boundaries are obtained. In fact, for arbitrary $\alpha > \frac{N-1}{N}$, the existence of  solutions  $(E_1, \dots, E_k) $ for both problems \eqref{f:pb1.b}-\eqref{f:pb2.b} is quite immediate,
and for convenience of the reader it  will be briefly discussed in Section \ref{sec:proofs} below. On the other hand, the analysis of their qualitative properties may require some more attention (in particular for the {\it maximum} problem \eqref{f:pb1.b}), but it is not our purpose to discuss here regularity issues.

Let us also mention that
the asymptotics  as $k \ra +\infty$ of the energies in \eqref{f:pb1.b}-\eqref{f:pb2.b}   has received a lot of attention in some recent works focused on  the honeycomb conjecture. Actually this celebrated conjecture, which was proved by Hales \cite{Hales} for cells of equal area minimizing the total perimeter, was already formulated in the context of the classical Cheeger constant in \cite{Car17}, in relationship with a conjecture by  Caffarelli and Lin involving optimal partitions of spectral type \cite{CaffLin}.
A proof of the honeycomb conjecture for the $\alpha$-Cheeger constant  when $\alpha = 1$ and $\alpha = 2$  has been  obtained very recently, respectively in \cite{bfvv17} and in \cite{bf17R}, under the restriction the  the admissible clusters are made by convex cells.

In this paper, we are interested in the study of $\alpha$-Cheeger clusters. Our aim is twofold. A first purpose is to get qualitative results describing the behaviour of $\alpha$-Cheeger clusters, in the limit when  $\alpha \ra \Big (\frac{N-1}{N}\Big )_+$ or $\alpha \ra +\infty$. A second purpose is to give an efficient phase field numerical approach for the computation of optimal $\alpha$-Cheeger clusters and, as a consequence of their asymptotic behaviour, of  optimal packings of balls in arbitrary boxes $D$.
The theoretical results are presented and discussed in the two subsections hereafter. Proofs and numerical results are then given respectively in Sections \ref{sec:proofs} and \ref{sec:num}.

\medskip

\subsection{Limiting behaviour of $\alpha$-Cheeger clusters} 
Our main asymptotical results show that, as $\alpha \ra \Big (\frac{N-1}{N}\Big )_+$, solutions to 
problems \eqref{f:pb1.b} and \eqref{f:pb2.b} converge to solutions to two different optimal packing problems for balls. More precisely: solutions of problem \eqref{f:pb1.b} converge to  a solution of the classical packing problem, which consists in finding $k$ mutually disjoint equal balls with maximal radius in $D$; solutions of problem \eqref{f:pb2.b} converge to  a solution of a more peculiar packing problem, which consists in finding $k$ mutually disjoint balls in $D$ maximizing the product of their volumes. 
The statements read as follows, where the $L ^1(D, \R^k)$-convergence has to be meant as the convergence of the characteristic functions $\{1_{\Om_1^\alpha}, \dots, 1_{\Om_k^\alpha}\}$ of solutions of \eqref{f:pb1} or \eqref{f:pb2} to the characteristic functions of balls $\{1_{B (x_1, r_1)}, \dots, 1_{B (x_k, r_k)}\}$ solving \eqref{op1} or \eqref{op2}.

\begin{theorem}\label{t:approx1}
As $\alpha  \rightarrow \Big (\frac{N-1}{N}\Big )_+$, up to subsequences, a solution to
\begin{equation}\label{f:pb1}\min  \Big \{\max_{i=1, \dots, k }   \frac{{\mathcal H} ^ {N-1}(\partial^* \Omega_i)} {|\Om_i| ^ {\alpha} }\   :\   \{\Omega _i\} \in \mathcal P _k (D)  \Big \}
\end{equation}
converges in $L ^1(D, \R^k)$  to a family of balls  solving the following optimal packing problem
\smallskip
\begin{equation}\label{op1}
\max \Big \{ r \ :\ \exists \{x_i , r_i \geq r\}_{i=1, \dots, k} \, , \ B (x_i, r_i) \subset D\, , \ B (x_i, r_i) \cap B (x_j , r_j)  = \emptyset \Big \}\,.
\end{equation}
\end{theorem}
\bigskip

\begin{theorem}\label{t:approx2}
As $\alpha  \rightarrow \Big (\frac{N-1}{N}\Big )_+$, up to subsequences, a solution to
\begin{equation}\label{f:pb2}\min  \Big \{\sum_{i=1} ^ k    \frac{{\mathcal H} ^ {N-1}(\partial^* \Omega_i)} {|\Om_i| ^ {\alpha} }\   :\   \{\Omega _i\} \in \mathcal P _k (D)  \Big \}\end{equation}
converges in $L ^1(D, \R^k)$  to a family of balls  solving the following optimal packing problem
\smallskip
\begin{equation}\label{op2}
\max \Big \{ \prod_{i=1} ^k r_i  \ :\ \exists \{x_i , r_i \}_{i=1, \dots, k} \, , \ B (x_i, r_i) \subset D\, , \ B (x_i, r_i) \cap B (x_j , r_j)  = \emptyset \Big \}\,.
\end{equation}
\end{theorem}
\bigskip

When $\alpha \ra +\infty$, we are not able to give a complete description of the asymptotic behaviour of $\alpha$-Cheeger clusters. Nevertheless we observe that, in case of problem \eqref{f:pb1},  such behaviour seems to be related to an optimal partition problem into cells of equal volume which minimize the product of their perimeters. This question is naturally linked to  the possible validity of a stronger version of the classical honeycomb conjecture solved by Hales. 
The picture is detailed in the statement given hereafter, along with a conjecture. For every $k \in \N$, we call $k$-cell a connected region obtained as the union of $k$ unit area regular hexagons taken from the hexagonal tiling of $\R ^2$. 

\begin{proposition}\label{l:infty}
As $\alpha \to + \infty$, up to subsequences a solution to problem \eqref{f:pb1} converges in $L ^ 1(D, \R^k)$ to a partition of $D$ into $k$ mutually disjoint subsets of equal measure. Moreover, if $N=2$, $D$ is a $k$-cell, and Conjecture \ref{Halesp} below holds true, this partition is the hexagonal one. 
\end{proposition}

\begin{conjecture}\label{Halesp}
Let $D$ be a $k$-cell. Then we expect that Hales' celebrated result \cite{Hales} can be strengthened into the following product-version of the honeycomb conjecture:
$$\min  \Big \{ \prod _{i=1} ^ k \mathcal H ^ 1(\partial^* \Omega _i) \ :\  \{ \Omega _i \} \in \mathcal P _k (D) \, ,  \  |\Omega _i |= 1 \Big \} = \big ( \mathcal H ^ 1( \partial H) \big ) ^ k  \,.$$

\end{conjecture}

\begin{figure}[!ht]
\includegraphics[width=0.25\textwidth]{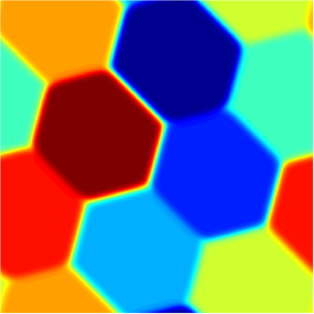}~
\includegraphics[width=0.25\textwidth]{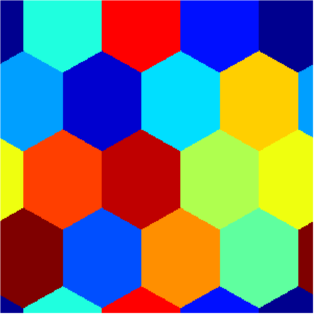}~
\includegraphics[width=0.25\textwidth]{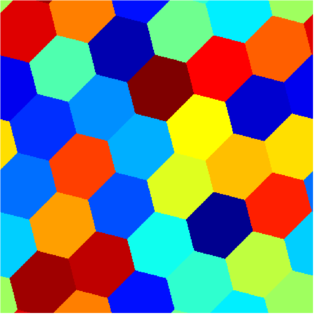}
\caption{Illustration of Conjecture \ref{Halesp} for $8,16,32$ cells in the torus}
\label{perim_prod}
\end{figure}

The proof of this new conjecture  seems to be non-trivial, as it would require a product form of the hexagonal isoperimetric inequality in the spirit of Hales. It is not a purpose of this paper to discuss such issue in detail.
We limit ourselves to support the validity of the conjecture by performing a few simulations using algorithms similar to those introduced in \cite{Oudet11}, the difference coming 
from the presence of the logarithmic function. In fact, in order to avoid numerical instabilities, we minimize the logarithm of the product of perimeters:
$$\min  \Big \{ \sum _{i=1} ^ k \log \big( \mathcal H ^ 1(\partial^* \Omega _i)  \big ) \ :\  \{ \Omega _i \} \in \mathcal P _k (D) \, ,  \  |\Omega _i |= 1 \Big \}.$$
The numerical results presented in Figure \ref{perim_prod} are performed for $8, 16$ and $32$ cells in the torus (the choice of the torus instead of a $k$-cell is based on the observation that,  if Conjecture \ref{Halesp} would fail on some $k$-cell which tiles the plane, then it will also fail on the torus). 

Concerning the asymptotic behaviour of $\alpha$-Cheeger clusters for problem \eqref{f:pb2} 
when $\alpha \ra +\infty$, the same arguments used to obtain  Proposition \ref{l:infty} can be used to show that, up to subsequences, a solution to such problem converges as well in $L ^ 1$ to a partition of $D$ into $k$ mutually disjoint subsets of equal measure. Anyway in this case we have no conjecture about the optimization problem solved by  this limit configuration.


\medskip
\subsection{Phase field approach for computing  $\alpha$-Cheeger clusters and optimal packings} There exist various works in the literature dealing with the computation of the Cheeger sets of a given domain $D$ 
(which corresponds to take $\alpha = 1$ and $k = 1$), using techniques related to convex optimization as, for instance, \cite{cfm09}. In \cite{CCP09} the authors exploit a projection algorithm in order to solve a more general weighted Cheeger problem, while in \cite{LRO05} a convex hull method is used for the optimization in the class of convex bodies. In \cite{KLR06} a characterization of the Cheeger set for convex domains is given, which can easily provide an algorithm for the computation of the Cheeger set for convex polygons. To our knowledge, no algorithm was implemented for the numerical computation of Cheeger clusters. An algorithm for computing numerically partitions of minimal length 
was introduced in \cite{Oudet11},
based on the approximation of perimeter by $\Gamma$-convergence using the classical Modica-Mortola theorem \cite{modica-mortola}. 

Our approach is as well of $\G$-convergence type and relies on the Modica-Mortola theorem for the approximation of the perimeter term, at the numerator in the $\alpha$-Cheeger ratio. Nevertheless, two difficulties arise: on one hand, even for the computation of the $\alpha$-Cheeger constant of a set (corresponding to the case $k=1$), one has to handle the approximation of the measure term in the denominator. This is done by considering a sufficiently high $L^p$-norm of the phase field function, as the more natural $L^1$ or $L^2$ norms are not suitable. Precisely, the measure of the set is efficiently approximated by the $L^\frac{2N}{N-1}$-norm, this choice being available for every $\alpha$. On the other hand, one has to cope with the collective behaviour of the different phases and the empty regions,  this being handled in the spirit of the partitioning algorithm introduced in \cite{BBO09}, with the novelty of the presence of the empty regions.  The outcoming $\Gamma$-convergence result reads: 

\begin{theorem} \label{t:mm}
For any fixed $\alpha >\frac{N-1}{N} $ and $p>1$ , the functional
$$F(\Om_1, \dots, \Omega_k):= 
 \sum_{i=1}^k \Big (\frac{{\mathcal H}^{N-1}(\partial^* \Omega_i)}{|\Omega_i|^\alpha} \Big ) ^p $$
is the $\Gamma$-limit as $\vps \ra 0$ in $L^1(D, \R^k)$ of 
the sequence 
$$F_\vps (u_1, \dots, u_k):=  \sum_{i=1}^k \Big ( \frac{\vps \int_D|\nabla u_i|^2 dx+  \frac{9} {\vps}\int _D u_i^2 (1-u_i)^2}{\Big (\int_D |u_i|^\frac{2N}{N-1}dx \Big )^\alpha} \Big ) ^ p$$\
\hskip 3.5 cm if $u_i \in H_0^1(D)$, $ u_i \ge 0$,  $\sum _{i=1} ^k u _i \leq 1$,  and $+\infty $ if not\,.
\end{theorem}

In the above statement, the separation of phases is related to the constraints  $u_i \ge 0$, $\sum _{i=1}^k u_i \le 1$, and to the fact that each phase will naturally converge to a bang-bang configuration. 
In order to avoid handling the constraint $\sum _{i=1}^k u_i \le 1$, we use a penalization approach inspired by the paper of Caffarelli and Lin \cite{CaffLin}.  Namely, we slightly change our functional by introducing the penalty term  $\frac 1\vps \sum_{1 \le i < j \le k} \int_D u_i^2u_j^2dx$. 
This term enhances the separation of phases during the computational process and avoids (at least the)
 local minima consisting on flat  functions at an intermediate level between $0$ and $1$. 
At the same time, the presence of this term, together with the bounds  $0\le u_i \le 1$, still ensure that the $\Gamma$-convergence holds. 
Thus we obtain the following corollary, which is 
the main practical tool for the numerical implementation our phase field approach:
\begin{corollary} \label{t:mm1}
For any fixed $\alpha >\frac{N-1}{N} $ and $p>1$, the functional
$$F(\Om_1, \dots, \Omega_k):= 
 \sum_{i=1}^k \Big (\frac{{\mathcal H}^{N-1}(\partial^* \Omega_i)}{|\Omega_i|^\alpha} \Big ) ^p $$
is the $\Gamma$-limit as $\vps \ra 0$ in $L^1(D, \R^N)$ of 
the sequence 
\begin{equation} 
F_\vps (u_1, \dots, u_k):=  \sum_{i=1}^k \Big ( \frac{\vps \int_D|\nabla u_i|^2 dx+  \frac{9} {\vps}\int _D u_i^2 (1-u_i)^2}{\Big (\int_D |u_i|^\frac{2N}{N-1}dx \Big )^\alpha} \Big ) ^ p+ \frac 1\vps \sum_{1 \le i < j \le k} \int_D u_i^2u_j^2dx
\label{optim_func}
\end{equation}
\
\hskip 3.5 cm if $u_i \in H_0^1(D)$,   $0\leq u_i \leq 1$  and $+\infty $ if not\,.
\end{corollary}

Our algorithm is very efficient in computing clusters or just Cheeger sets in arbitrary boxes $D$, in two or three space dimensions. In Section \ref{sec:num}, we give several numerical experiments with different values of $\alpha$ and $k$,  in $\R ^2$ and in $\R ^3$. In particular, for $\alpha$ close to the critical value $\Big (\frac{N-1}{N}\Big )_+$, we recover some classical results for optimal circle/sphere packings, underlining the interest of Theorems \ref{t:approx1} and \ref{t:approx2}. The literature related to circle packings is quite vast, and efficient algorithms are known to work for thousands of cells. However, they rely on heuristic combinatorial and geometric ideas in order to search for the optimal configuration (see for instance \cite{cpack10, cpackLB11}). Our approach is a global one: we initialize each cell with a random density function and we perform a direct gradient based optimization to reach the final configuration. This approach works rather well in the cases we considered. We underline the fact that the search of optimal circle packings is a non-smooth problem. The desired configuration corresponds to maximizing the minimal distance between centers while minimizing the size of the container. Any standard gradient descent algorithm will get stuck as soon as the pairwise distances between centers contain repeated minimal values. The method we present can consistently get near an optimal candidate, or even produce a good starting point for other algorithms.
\medskip

%
%

\smallskip

\section{Proofs}\label{sec:proofs}

Throughout this section, we set $\omega _N :={\mathcal H ^ {N-1} (\partial B)}/{|B| ^ { \frac{N-1}{N}}}$, being $B$ a ball in $\R ^N$. 

\medskip
 {\bf Existence of $\alpha$-Cheeger clusters}.
Before giving the proofs of our results, for the benefit of the reader we briefly discuss the existence of solutions to problems \eqref{f:pb1.b} and \eqref{f:pb2.b} for a  fixed $\alpha > \frac{N-1}{N}$.  This question relies on classical compactness results in $BV$-spaces and, for $\alpha=1$ and problem \eqref{f:pb2.b}, has been extensively discussed in \cite{Car17}. For shortness, we only deal with problem \eqref{f:pb1.b} (the arguments being the same for problem \eqref{f:pb2.b}). Assume that $\{\Om_1^n, \dots, \Om_k^n\}$ are elements of $\mathcal P_k(D)$ such that
$$\max_{i=1, \dots, k }   \frac{{\mathcal H} ^ {N-1}(\partial^* \Omega^n_i)} {|\Om^n_i| ^ {\alpha} } \longrightarrow \inf   \Big \{\max_{i=1, \dots, k }   \frac{{\mathcal H} ^ {N-1}(\partial^* \Omega_i)} {|\Om_i| ^ {\alpha} }\   :\   \{\Omega _i\} \in \mathcal P _k (D)  \Big \}\;\; \text{ as } \;\; n \rightarrow +\infty.$$
One has first to observe that 
\begin{equation}\label{bbf003}
\liminf_{n \rightarrow +\infty} |\Om_i^n| >0 \qquad \forall i=1, \dots, k\,.
\end{equation}
Indeed, we can assume that there exists some $M>0$ such that $\max_{i=1, \dots, k }   \frac{{\mathcal H} ^ {N-1}(\partial^* \Omega^n_i)} {|\Om^n_i| ^ {\alpha} }\le M$. Then, using the isoperimetric inequality on each cell $\Omega^n_i$ we get
$$\omega_N \frac{1}{|\Om^n_i|^{\alpha - \frac{N-1}{N}}}\le \frac{{\mathcal H} ^ {N-1}(\partial^* \Omega^n_i)} {|\Om^n_i| ^ {\alpha} }\le M,$$
hence, \eqref{bbf003} is true. Exploiting also the fact that 
$${\mathcal H} ^ {N-1}(\partial^* \Omega^n_i) \le M |D|^k,$$
we get that the sequence $(1_{\Om_i^n})_n$ is bounded in $BV (D)$. By standard compactness results, we can assume that,  up to a subsequence,  it converges strongly in $L^1(D) $ to some limit which can be written as $1_{\Om_i}$. For $i\not =j$, the sets $\Om_i$ and $\Om_j$ are disjoint and, in view of \eqref{bbf003}, they have strictly positive measure. Moreover,  by lower semicontinuity of the perimeter, we have 
$$\liminf_{n \rightarrow +\infty} {\mathcal H} ^ {N-1}(\partial^* \Omega^n_i) \ge  {\mathcal H} ^ {N-1}(\partial^* \Omega_i),$$
hence  $\{\Om_1, \dots, \Om_k\}$ is a solution to  problem \eqref{f:pb1.b}.

\bigskip

{\bf Proof of Theorem \ref{t:approx1}}. 
We write $\alpha = \frac{N-1}{N} + \delta$ for some $\delta >0$, and we consider the auxiliary problems 
$$\begin{array}{ll} M_{k , \delta} (D)& \displaystyle:=\min  \Big \{\max_{i=1, \dots, k }   \left ( \frac{{\mathcal H} ^ {N-1}(\partial^* \Omega_i)} {\omega _N |\Om_i| ^ { \frac{N-1}{N} + \delta} } \right ) ^ {\frac{1}{\delta}}\   :\   \{\Omega _i\} \in \mathcal P _k (D)  \Big \}
\\ \noalign{\medskip} 
& \displaystyle =\min  \Big \{\max_{i=1, \dots, k }   \left ( \frac{{\mathcal H} ^ {N-1}(\partial^* \Omega_i)} {\omega _N |\Om_i| ^ { \frac{N-1}{N}} } \right ) ^ {\frac{1}{\delta}}\!\! \frac{1}{|\Omega _i |}    :\   \{\Omega _i\} \in \mathcal P _k (D)  \Big \}\,, 
 \end{array}$$
which have the same optimal clusters as our  initial problems \eqref{f:pb1}. 

If  $D$ contains $k$ mutually disjoint balls of radius $r_D$, we have 
\begin{equation}\label{firstbound}
M _{k , \delta} (D) \leq \frac{1}{|B _{r_D}|} \,.
\end{equation}
Then, if  $\{ \Omega ^ 1 _\delta, \dots \Omega ^ k _ \delta \}$ is an optimal cluster for $M _k ^ \delta (D)$, we have
\begin{equation}\label{f:estinf} 
\left ( \frac{{\mathcal H} ^ {N-1}(\partial^* \Omega^ \delta _i)} {\omega _N |\Om ^ \delta _i| ^ { \frac{N-1}{N}} } \right ) ^ {\frac{1}{\delta}}\!\! \frac{1}{|\Omega _i^ \delta |} \leq \frac{1}{|B_{r_D}|} \qquad \forall i = 1, \dots, k \,.
\end{equation}   
Using \eqref{f:estinf} and  the isoperimetric inequality, we get 
\begin{equation}\label{f:boundarea}
{|\Omega _i^ \delta |} \geq {|B_{r_D}|} \qquad \forall i = 1, \dots, k\,.
 \end{equation} 
Moreover, 
using \eqref{f:estinf} and the upper bound $|\Om ^ \delta _i| \leq |D|$, we get 
\begin{equation}\label{f:boundper}
{\mathcal H} ^ {N-1}(\partial^* \Omega^ \delta _i) \leq \omega _N |D| ^ { \frac{N-1}{N}} \Big ( \frac{|D|}{|B_{r_D}|}  \Big )^ {\delta}\qquad \forall i = 1, \dots, k\,. 
\end{equation}
We deduce that a (not relabeled) subsequence of $\{ \Omega _ 1^\delta, \dots \Omega ^  \delta_k \}$ converges in $L^1(D, \R^k)$ to a limit cluster, 
that we denote by $\{ \Omega _1, \dots, \Omega _ k \}$, satisfying $|\Omega _i| \geq |B_{r_D}|$. Let us show first that all the sets $\Omega _i$'s are balls, and then that they solve the optimal packing problem \eqref{op1}. 

By using the definitions of $M _{k, \delta} (D)$  and $r _D$,  we have the following estimate:
$$\frac{{\mathcal H} ^ {N-1}(\partial^* \Omega_i ^ \delta )} { |\Om_i ^ \delta | ^ { \frac{N-1}{N} + \delta} } \leq \omega _N (M _{k, \delta} (D) )^ \delta
\leq\frac {\omega _N}{|B_{r_D}| ^ \delta}\qquad \forall i = 1, \dots, k\,.
 $$
Then, passing to the liminf as $\delta \to 0$, we obtain 
$$\frac{{\mathcal H} ^ {N-1}(\partial^* \Omega_i  )} { |\Om_i  | ^ { \frac{N-1}{N}} } \leq \omega _N \qquad \forall i = 1, \dots, k \,,$$ 
which implies that the sets $\Omega _i$'s are balls. 

Let us show that they solve the optimal packing problem \eqref{op1}. Denote by $r_1, \dots r _k$ the radii of $\Omega _1, \dots \Omega _k$, set $r_{\min} := \min \{ r_1, \dots , r _k \}$ and let $r _*$ denote the maximal radius  in the optimal packing problem \eqref{op1}. 
Clearly, by the definition of $r_*$, it holds $r _* \geq r _{\min}$. 
 On the other hand, since $D$ contains $k$ mutually disjoint balls of radius $r_*$, by definition of $M _{k, \delta} ( D)$ we have 
$$M _{k, \delta} (D) \leq \frac{1}{ |B _{r _* }|}\,.$$
Hence the inequalities \eqref{f:estinf} and \eqref{f:boundarea} are in force with $r_*$ in place of $r_D$ and therefore, passing to the limit as $\delta \to 0$, we obtain
$$|\Omega _i | \geq |B _{r_*}|\qquad \forall i = 1, \dots, k \,. $$
Then $r_{\min } \geq r_* $ and we conclude that $r _{\min } = r _* $. \qed

\bigskip\bigskip
{\bf Proof of Theorem \ref{t:approx2}}. 
We write $\alpha = \frac{N-1}{N} + \delta$  for some $\delta >0$, and we set 
$$
m_{k , \delta} (D):=\min  \Big \{\sum_{i=1}^  k     \frac{{\mathcal H} ^ {N-1}(\partial^* \Omega_i)} {\omega _N |\Om_i| ^ { \frac{N-1}{N} + \delta} }\   :\   \{\Omega _i\} \in \mathcal P _k (D)  \Big \} \,,$$
(where we have introduced just  for convenience the constant $\omega _N$ with respect to our  initial problems \eqref{f:pb2}). 
 
Let $(\Omega _1 ^ \d, \dots, \Omega ^\d _k )$ be an optimal cluster for
$m_{k , \delta} (D)$. From the inclusion $\Omega _i ^ \d \subset D$, we know that
\begin{equation}\label{f:limsup}
\limsup _{\d \to 0} |\Omega _i ^ \d | ^ \d \leq \limsup _{\d \to 0} |D  | ^ \d   = 1  \qquad \forall i = 1, \dots, k \,. 
\end{equation}
We claim that 
\begin{equation}\label{f:lim}
 \lim _{\d \to 0} |\Omega _i ^ \d | ^ \d   = 1 \qquad \forall i = 1, \dots, k \,.
\end{equation}
Indeed, if  $D$ contains $k$ mutually disjoint balls of radius $r_D$, we have 
\begin{equation}\label{firstbound2}
\sum_{i=1}^  k     \frac{1}{ |\Om_i^\d | ^ \d}  
\leq
\sum_{i=1}^  k     \frac{{\mathcal H} ^ {N-1}(\partial^* \Omega^\d _i)} {\omega _N |\Om_i^\d | ^ { \frac{N-1}{N} } } \frac{1}{ |\Om_i^\d | ^ \d} = m _{k , \delta} (D) \leq \frac{k}{|B_{r_D}|^\d } \,.
\end{equation}
In the limit as $\d \to 0$, up to passing to a subsequence such that $\{ |\Omega _i ^ \d | ^ \d\} $ converges, we find
$$ \sum _{i=1}^  k \frac{1}{ \lim _{\d \to 0} |\Omega _i ^ \d | ^ \d} \leq k \,.$$
In view of \eqref{f:limsup},   the above inequality implies that necessarily \eqref{f:lim} is satisfied. 

Now we observe that, by  the isoperimetric inequality, there holds
$$ \begin{array}{ll}
&  \displaystyle     \frac{{\mathcal H} ^ {N-1}(\partial^* \Omega^\d _1)} {\omega _N |\Om_1^\d | ^ { \frac{N-1}{N} } } \frac{1}{ |\Om_1^\d | ^ \d}  
+ \sum_{i=2}^  k        \frac{1}{ |\Om_i^\d | ^ \d}  
\leq
\\ \noalign{\medskip}
& \displaystyle \frac{{\mathcal H} ^ {N-1}(\partial^* \Omega^\d _1)} {\omega _N |\Om_1^\d | ^ { \frac{N-1}{N} } } \frac{1}{ |\Om_1^\d | ^ \d}  
+ \sum_{i=2}^  k       \frac{{\mathcal H} ^ {N-1}(\partial^* \Omega^\d _i)} {\omega _N |\Om_i^\d | ^ { \frac{N-1}{N} } } \frac{1}{ |\Om_i^\d | ^ \d} 
= m _{k , \delta} (D) \leq \frac{k}{|B_{r_D}|^\d }\,. \end{array}
 $$
In the limit as $\d \to 0_+$, thanks to \eqref{f:lim}, we deduce that 
\begin{equation}\label{f:uno}
\limsup _{\d \to 0} \frac{{\mathcal H} ^ {N-1}(\partial^* \Omega^\d _1)} {\omega _N |\Om_1^\d | ^ { \frac{N-1}{N} } } \leq 1 \,,
\end{equation} 
which implies in particular  
\begin{equation}\label{f:due}
\limsup _{\d \to 0} {{\mathcal H} ^ {N-1}(\partial^* \Omega^\d _1)}  \leq {\omega _N |D | ^ { \frac{N-1}{N} } }  \,.
\end{equation}
By repeating the same argument  for every $i = 2, \dots, k$, we obtain that
that the perimeters of $\Omega ^ \d _i$ are bounded from above for every $i = 1, \dots, k$. 

Next, let us show that the measures 
 of $\Omega ^ \d _i$ are bounded from below for every $i = 1, \dots, k$.  
 We exploit once again the isoperimetric inequality to write 
 $$  \sum_{i=1}^  k  \exp \Big (\d \log      \frac{1}{ |\Om_i^\d | } \Big )  = \sum_{i=1}^  k       \frac{1}{ |\Om_i^\d | ^ \d} \leq  \sum_{i=1}^  k       \frac{{\mathcal H} ^ {N-1}(\partial^* \Omega^\d _i)} {\omega _N |\Om_i^\d | ^ { \frac{N-1}{N} } } \frac{1}{ |\Om_i^\d | ^ \d} 
= m _{k , \delta} (D) \leq \frac{k}{|B_{r_D}|^\d } \,.
 $$ 
 By using the elementary inequality $e ^ t \geq 1 + t$ for every  $ t \in \R$, we deduce that 
 $$k - \d  \log \Big ( \prod_{i =1} ^ k  |\Om_i^\d | \Big )  \leq\frac{k}{(\gamma _N r_D ^N) ^\d } \,, $$
or equivalently
 $$ \log \Big ( \prod_{i =1} ^ k  |\Om_i^\d | \Big )  \geq {k} \frac{1}{\d}\Big ( 1 - \frac{1}{(\gamma _N r_D^N) ^\d } \Big )\, ,    $$
 where $\gamma _N$ denotes the measure of the unit ball in $\R^N$. 

Since the right hand side of the above inequality remains bounded as a function of $\d$ in a neighbourhood of $\delta = 0$, 
and since logarithm is monotone increasing, we infer that the product of the measures of $\Om ^ \d_i$'s remains bounded from below. 
Since each of these measures is bounded from above (by the measure of $D$), we conclude that 
$|\Omega _i ^ \d|$ is bounded below for every $i = 1, \dots, k$. 

So far, we have obtained that,  for every $i = 1, \dots, k$, the sets $\Omega _ i ^\d$ have perimeters bounded from above and measures bounded from below {\it uniformly} as $\delta\ra 0_+$; hence, up to a (not relabeled) subsequence, they admit  a $k$-uple  of (non trivial) limit sets, and from inequalities \eqref{f:uno} we see that these limit sets are necessarily balls. Let us denote them by $(B^0 _1, \dots, B^0 _k)$, and let us show that they solve the optimal packing problem \eqref{op2}. 

To that aim, we are going to estimate from above and from below the quotient 
$$\frac{m _{k, \d} (D) -k }{\d} \,.$$ 
On one hand, if $( B ^*_1, \dots, B ^*_k)$ are balls which solve the optimal packing problem \eqref{op2}, taking them as a test cluster in the definition of $m _{k, \d} (D)$, we get the upper  bound 
$$
\frac{m _{k, \d} (D) -k }{\d} \leq  \frac{1}{\d} \sum_{i=1}^  k     \Big (    \frac{1}{ |B ^*_i | ^ \d} -1\Big )  = 
\sum_{i=1}^  k  \frac{1}{\d} \Big ( \exp \Big (\d \log      \frac{1}{ |B ^* _i | } \Big )  - 1 \Big)\,.
 $$Hence, 
\begin{equation}\label{UB}
\limsup _{\d \to 0} \frac{m _{k, \d} (D) -k }{\d} \leq \log \Big ( \prod_{i=1} ^ k\frac{1}{ |B ^*_i | }\Big ) \,. 
\end{equation}
On the other hand, by applying as usual the isoperimetric inequality, we have
$$\frac{m _{k, \d} (D) -k }{\d} \geq \frac{1}{\d} \sum_{i=1}^  k     \Big (    \frac{1}{ |\Omega_i ^\d | ^ \d} -1\Big )\,.$$ 
Now we exploit the fact that $ |\Omega_i ^\d | \to |B _i ^0|$ as $\d \to 0_+$, so that,  for every fixed $\eta>0$, we can find $\overline \d = \overline \d (\eta)$ such that 
$$|\Omega_i ^\d | \leq |(1+ \eta) B_i^0 | \qquad \forall \d \leq \overline \d (\eta)\, , \ \forall i = 1, \dots, k\,.$$   
Therefore, 
$$\frac{m _{k, \d} (D) -k }{\d} \geq \frac{1}{\d} \sum_{i=1}^  k     \Big (    \frac{1}{|(1+ \eta) B_i^0 | ^ \d} -1\Big )\,.$$ 
Then, by arguing as above and using the arbitrariness of $\eta$, we arrive at
\begin{equation}\label{LB}
\limsup _{\d \to 0} \frac{m _{k, \d} (D) -k }{\d} \leq \log \Big ( \prod_{i=1} ^ k\frac{1}{ |B ^0_i | }\Big ) \,. 
\end{equation} 
Eventually, combining \eqref{UB} and \eqref{LB}, we obtain 
$$ \prod_{i=1} ^ k\frac{1}{ |B ^0_i | } \geq  \prod_{i=1} ^ k\frac{1}{ |B ^*_i | }\,.$$
By the definition of $B ^ * _i$, the above inequality holds necessarily with equality sign, which amounts to say that the limit balls $B^0_i$ maximize the product of volumes  of a family of $k$ mutually disjoint balls contained into $D$, namely they solve problem \eqref{op2}. \qed

\bigskip\bigskip
{\bf Proof of Proposition \ref{l:infty}}. 
Let $(\Omega ^ \alpha _1, \dots , \Omega ^ \alpha _k)$ by a solution to problem \eqref{f:pb1}, and assume $\alpha \rightarrow +\infty$. 
Assuming without loss of generality that $|D| = k$, we consider a $k$-uple $(\widehat \Omega  _1, \dots, \widehat \Omega _k) \in \mathcal P _k ( D)$ of  sets with finite perimeter, such that $|\widehat \Omega  _i | = 1$ for every $i = 1, \dots, k$.  Then we have
\begin{equation}\label{f:2334}
\frac{{\mathcal H} ^ {N-1}(\partial^* \Omega^\alpha _i)} {|\Om_i ^\alpha| ^ {\alpha} } \leq C:= \max _{i = 1, \dots, k} |\partial^* \widehat \Omega  _i| \qquad \forall i = 1, \dots, k\,.
\end{equation}
By applying the isoperimetric inequality, we obtain 
$$\frac{ \omega _N |  \Omega^\alpha _i| ^ {\frac{ N-1}{N}}} {|\Om_i ^\alpha| ^ {\alpha} } \leq C \qquad \forall i = 1, \dots, k\,,$$ 
or equivalently 
$$|\Om_i ^\alpha| \geq \Big ( \frac{\omega _N}{C} \Big ) ^ { \frac{1}{\alpha - \frac{N-1}{N}}} 
 \qquad \forall i = 1, \dots, k\,.$$ 
Hence, 
\begin{equation}\label{f:2340}
\liminf_{\alpha \to + \infty} |\Om_i ^\alpha| \geq 1 \qquad \forall i = 1, \dots, k\,.
\end{equation}
Moreover, from inequality \eqref{f:2334}, we know that
$$\sum _{i = 1} ^ k \Big (\frac{\mathcal H ^ {N-1} (\partial^* \Omega ^ \alpha _i )}{C} \Big ) ^ {\frac{1}{\alpha}} \leq 
\sum _{i = 1} ^ k   |  \Omega^\alpha _i| \leq k\,.$$ 
Thus, by using the elementary inequality $e ^ {t} \geq 1 + t$ for every $t \in \R$, we arrive at
$$\frac{1}{\alpha} \sum _{i = 1} ^ k \log  \Big (\frac{\mathcal H ^ {N-1} (\partial^* \Omega ^ \alpha _i )}{C} \Big ) \leq 0\,, $$ 
which means that
\begin{equation}\label{f:2348}
\prod_{i=1} ^ k \mathcal H ^ {N-1} (\partial^* \Omega ^ \alpha _i ) \leq C ^ k \,.
\end{equation}
We observe that,  by \eqref{f:2340},  $\mathcal H ^ {N-1} (\partial^* \Omega ^ \alpha _i )$ is bounded from below for every $i= 1, \dots, k$; then inequality \eqref{f:2348} ensures that $\mathcal H ^ {N-1} (\partial^* \Omega ^ \alpha _i )$ is also bounded from above for every $i= 1, \dots, k$. 
We conclude that $(\Omega ^ \alpha _1, \dots , \Omega ^ \alpha _k)$ admits a limit in $L ^1$, hereafter denoted by $(\Omega ^ \infty _1, \dots , \Omega ^ \infty _k)$, and in view of \eqref{f:2340} this limit turns out to be a partition of $D$.

Finally, if $N=2$ and $D$ is a $k$-cell, by taking the sets $\widehat \Omega _i$'s in \eqref{f:2334} equal to the $k$-copies of $H$ which compose the $k$-cell, inequality \eqref{f:2348} implies 
\begin{equation}\label{f:2349}
\prod_{i=1} ^ k \mathcal H ^ {1} (\partial^* \Omega ^ \infty _i ) \leq  \big ( \mathcal H ^ {1} (\partial H)  \big ) ^ k  \,.
\end{equation}
Therefore, if Conjecture \ref{Halesp} is satisfied, we conclude that
$$
\prod_{i=1} ^ k \mathcal H ^ {1} (\partial^* \Omega ^ \infty _i ) \leq \min  \Big \{ \prod _{i=1} ^ k \mathcal H ^ 1(\partial^* \Omega _i) \ :\  \{ \Omega _i \} \in \mathcal P _k (D) \, ,  \  |\Omega _i |= 1 \Big \}\,.
$$ 
Hence the above inequality holds with equality sign, and  the partition  $\{\Omega ^ \infty _i  \}$ is optimal for the minimization problem 
at the right hand side. \qed

\bigskip\bigskip
{\bf Proof of Theorem \ref{t:mm}}. We prove separately the so-called $\G$-liminf and $\G$-limsup inequalities (see the monograph \cite{DM} for an introduction to $\Gamma$-convergence). 

\noindent {\it $\G$-liminf inequality.} Let $u_i^\vps \in H_0^1(D)\sm \{0\}$, such that $(u_1^\vps,\dots, u_k^\vps)\sr{L^1(D, \R^k)}{\lra}(u_1, \dots, u_k)$ and
$$\limsup_{\vps \ra 0} F_\vps (u_1^\vps,\dots, u_k^\vps)<+\infty.$$
In a first step, we notice the existence of a constant $M$, such that for small $\vps$ and for every $i$ we have
$$\vps \int_D|\nabla u_i^\vps|^2 dx + \frac{9} {\vps}\int_D  (u_i^\vps)^2 (1-u_i^\vps)^2 dx \le M |D|^\alpha.$$
From the Modica-Mortola theorem, we get that $u_i \in BV(D, \{0,1\})$, hence $u_i =1_{\Om_i}$, where $\Om_i$ are pairwise disjoint and satisfy
$$\HN(\partial^*  \Om) \le \liminf_{\vps \ra 0} \vps \int_D|\nabla u_i^\vps|^2 dx + \frac{9} {\vps}\int_D  (u_i^\vps)^2 (1-u_i^\vps)^2 dx.$$
Let us prove that $\Om_i \not= \emptyset$. We know that
$$\vps \int_D|\nabla u_i^\vps|^2 dx + \frac{9} {\vps}\int_D  (u_i^\vps)^2 (1-u_i^\vps)^2 dx \le M \Big (\int_D (u_i^\vps)^\frac{2N}{N-1} dx\Big ) ^\alpha.$$
Setting $\delta:=  \frac{2N}{N-1}\alpha -2 >0$ and using the Cauchy-Schwartz inequality on the left hand side, we get
$$6 \int_D|\nabla u_i^\vps| u_i^\vps (1-u_i^\vps)dx \le M \|u_i^\vps\|^{2+\delta}_ {L^\frac{2N}{N-1}(D)},$$
or, by the chain rule,
$$6 \int_D\Big |\nabla \Big ( \frac{(u_i^\vps)^2}{2}- \frac{(u_i^\vps)^3}{3}\Big )\Big  | dx \le M \|u_i^\vps\|^{2+\delta}_ {L^\frac{2N}{N-1}(D)}.$$
Using the Sobolev inequality with a dimensional constant $S_N$, we get
$$6S_N \|\frac{(u_i^\vps)^2}{2}- \frac{(u_i^\vps)^3}{3}\| _{L^\frac{N}{N-1}(D)}\le M \|u_i^\vps\|^{2+\delta}_ {L^\frac{2N}{N-1}(D)}.$$
Since $0\le u_i^\vps \le 1$ we get $\frac{(u_i^\vps)^2}{2}- \frac{(u_i^\vps)^3}{3} \ge \frac{(u_i^\vps)^2}{6}$ so
$$ S_N \|u_i^\vps\|^2 _{L^\frac{2N}{N-1}(D)}\le M \|u_i^\vps\|^{2+\delta}_ {L^\frac{2N}{N-1}(D)},$$
or
$$ \frac{S_N}{M} \le \|u_i^\vps\|^{\delta}_ {L^\frac{2N}{N-1}(D)}.$$
Passing to the limit, we get $|\Om_i|\ge\frac{S_N}{M} $. Since  $\int_D (u_i^\vps)^\frac{2N}{N-1} dx \ra |\Om_i| \not =0$, the $\G$-liminf property occurs.
\medskip

\noindent {\it $\G$-limsup inequality.} For one single set, the Modica-Mortola theorem gives the procedure of constructing the recovering sequence. For partition problems we refer to the paper by Baldo \cite{BA90}, where the recovering sequence requires more attention because of the exact partition requirement which is obtained via the constraint $\forall x\in D, \; \sum_{i=1}^k u_i(x)=1$. Since in our problem we do not have a complete partition of the set $D$, we can give a direct proof as follows.

Let $(\Om_1, \dots, \Om_k)$ be pairwise disjoint measurable subsets of $D$ with finite perimeter. Let us fix a positive constant $\delta >0$. For a usual convolution kernel $(\rho_\eta)_\eta$, we have
\begin{equation}\label{bbf01} 
(\sum_{i=1}^k 1_{\Om_i} ) \ast \rho_\eta \le 1,
\end{equation}
and we claim that we can choose $\eta$ small enough such that for every $i=1, \dots, k$  we can find $t_i \in (\frac 12, 1)$ such that the set $A_i^\delta:= \{ 1_{\Om_i} \ast \rho_\eta >t_i\}$ is smooth and
$$\HN(\partial A_i^\delta) \le \HN(\partial^* \Om_i)+ \delta,$$
$$ |A_i^\delta \Delta \Om_i|\le \delta.$$
Indeed, from the Sard theorem, we know that almost all levels are smooth. From the strict convergence $\rho_\eta\ast 1_{\Om_i}\sr{BV}{\longrightarrow} 1_{\Om_i}$, we get that  
$$\int_0^1 |\nabla \rho_\eta\ast 1_{\Om_i}|dx \ra \HN(\partial ^*\Om_i),$$
and the co-area formula gives the existence of a level $t >\frac 12$ such that
$$\HN(\partial ^*\Om_i) \ge  \limsup_{\eta \rightarrow 0} \HN(\partial \{\rho_\eta\ast 1_{\Om_i}>t\}).$$

On the other hand, from \eqref{bbf01} and the choice of $t_i > \frac 12$, the sets $A_i^\delta$ are disjoint, at positive distance. Using now the Modica-Mortola theorem, we can find recovering sequences $u_i^\vps \sr{L^1(D)}{\ra} 1_{A_i^\delta}$ with 
$$\limsup_{\vps \ra 0} \vps \int_D|\nabla u_i^\vps|^2 dx + \frac{9} {\vps}\int_D  (u_i^\vps)^2 (1-u_i^\vps)^2 dx \le \HN(\partial A_i^\delta).$$
Choosing, a sequence $\delta _n\ra 0$, by a diagonal procedure we can find $\vps_n$ small enough such that the $\G$-limsup property holds, and moreover $u_i^{\vps_n}u_j^{\vps_n}=0$. \qed

\bigskip\bigskip

{\bf Proof of Corollary \ref{t:mm1}.} The proof is implicitly contained into the one  of  Theorem \ref{t:mm}. In particular, we point out that the choice of the recovering sequence for the $\G$-limsup property is also suitable for the penalized functionals.
\qed

\bigskip

\section{Numerical results}\label{sec:num}

In order to discretize the functional \eqref{optim_func}, we consider a rectangular box $D$ in $\Bbb{R}^2$ or $\Bbb{R}^3$ endowed with a finite differences uniform grid with $M$  discretization points along each axis direction. A function $u$ will be numerically represented by its values at the grid points. We use basic order $1$ centered finite differences in order to compute the gradient terms $|\nabla u_i|$, and basic quadrature formulas to compute all integrals. Similar approaches were already used in \cite{Oudet11} and \cite{BO16}. Moreover, in \cite[Section 4]{BO16} detailed expressions of the gradient with respect to each of the grid point variables of \eqref{optim_func} are given for the Modica-Mortola term. The other integrals, like the denominator of $F_\varepsilon$ and the penalization terms, are approximated by their arithmetic mean along the grid. The optimization is done using a LBFGS quasi-Newton method implemented in Matlab \cite{lbfgs}. The algorithm uses information on a number of previous gradients, $5$ in our computations, in order to build an approximation of the Hessian. In addition to being more efficient in avoiding local minima than a simple gradient descent algorithm, the LBFGS algorithm recalled above allows us to impose pointwise bounds for every variable, fact which is important in our approach. This also motivated us to use a penalized approach, rather than a projected gradient approach. The optimization procedure is presented in Algorithm \ref{alg:opt}.

\renewcommand{\labelitemi}{$\cdot$}
\begin{algorithm}
\caption{Optimization procedure}
\label{alg:opt}
\begin{algorithmic}[1]
\Require 
\begin{itemize}
\item $M$: initial discretization parameter ($M=20$ in 2D, $M=10$ in 3D)
\item $k$: number of cells
\item $d$: dimension
\item $maxit$: maximal number of iterations
\item $\varepsilon \in [1/M,4/M]$: Modica-Mortola parameter
\item $tol$: stopping criterion
\item $n$: number of refinements chosen so that the final resolution is as good as we want (in our case $2^{n-1}M>300$ in 2D and $2^{n-1}M>100$ in 3D)
\end{itemize}
\State Initialize densities $u_1,..,u_k$ as $k$ random vectors of $M^d$ elements
\State step $= 1$
\Repeat
  \If{step$>1$} 
    \State $M\gets 2M$
    \State interpolate linearly the previous densities $(u_i)_{i=1}^k$ on the new grid
  \EndIf 
  \State Run the LBFGS optimization procedure \cite{lbfgs} with the following parameters 
  \begin{itemize}
  \item starting point $(u_i)_{i=1}^k$ (random for first step, interpolated from the previous optimization result for the next steps)
  \item tolerances and maximum number of iterations: $tol=10^{-8},\ maxit=10000$
  \item function to optimize: \eqref{optim_func}. Value and gradient are computed at each iteration
  \item pointwise upper and lower bounds $0\leq u_i \leq 1$.
  \end{itemize} 
  \State The previous algorithm returns the optimized densities $(u_i)_{i=1}^n$

\State step = step$+1$ (go to next step)
\Until step$>n$

\Return the $k$ density functions
\end{algorithmic}
\end{algorithm}

Since we perform an optimization of a non-convex functional, there is no guarantee of convergence to a global minimum. In order to avoid local minima, we choose as initialization some random densities for each of the functions $u_i$. Moreover, in order to avoid the rapid convergence to a characteristic function and to diminish the effects of the non-convex potential, we choose the parameter $\varepsilon$ equal to $1/M$, where $M$ is the number of discretization points along each axis direction. The parameter $\varepsilon$ dictates the width of the interface where the functions $u_i$ go from $0$ to $1$. Sometimes it is useful to consider larger values of $\varepsilon$, $2/M$ or $4/M$, to allow the cells to move more freely. 

In order to reduce the number of iterations needed to reach an optimum, we propose a grid refinement procedure as already noted in \cite{BBO09} and \cite{BO16}. We perform an initial optimization on a grid of rather low size with $M \in [20,50]$. Then we interpolate the result on a finer grid, usually doubling the discretization parameter, and we continue the optimization on this refined grid. We continue until we reach the desired level of accuracy. In 2D we can easily perform computations on grids of size up to $400\times 400$ while in 3D we go up to $100\times 100 \times 100$.  Our current algorithm works well if the number of cells is not too large. When considering many phases, computations become more costly and memory costs become larger. It is possible to use techniques from \cite{BogLS17} in order to address these issues. The idea is that, when dealing with partitioning or multiphase problems in the phase-field context, one could restrict the computation in a neighborhood of the significant part of each cell, for example $\{u_i>\delta\}$, for a given threshold $\delta>0$. This could significantly reduce the computational and memory cost of the computations and could allow the use of the algorithm for many cells. However, this goes beyond the scope of this article.

Even if we choose to work on finite differences grid we may still compute $\alpha$-Cheeger sets and $\alpha$-Cheeger clusters corresponding to non-rectangular domains. For a general domain $D$ we consider a rectangular box $D' \supset D$ on which we construct the finite differences grid. We set all functions involved in the computations to be equal to zero on grid points outside $D$ and set the gradient to be equal to zero on the same points lying outside $D$. In this way the optimization is made only on points inside the desired domain $D$. 

It is also possible to use a finite element framework in order to minimize \eqref{optim_func} on general domains. In \cite[Section 3]{BO16} you can find a detailed presentation of such a finite element framework in the context of Modica-Mortola functionals. Once the mass and rigity matrices for the Lagrange P1 finite elements are obtained, all functionals needed in our computations can be expressed as vector matrix products.  

Now we exemplify the use of the proposed algorithm for computing $\alpha$-Cheeger sets, $\alpha$-Cheeger clusters and optimal packings. We make available an implementation of the algorithm described above which can be found online at the following link: \href{https://github.com/bbogo/Cheeger_patch}{\nolinkurl{https://github.com/bbogo/Cheeger_patch}}. This implementation uses the finite element framework for the optimization of \eqref{optim_func}.  As detailed below it is also possible, for convex domains, to compare the Cheeger sets found by minimizing \eqref{optim_func} with the exact Cheeger sets obtained by using the representation formula provided by Kawohl and Lachand-Robert in \cite{KLR06}. 
\medskip

\begin{figure}[ht]
\includegraphics[width=0.2\textwidth]{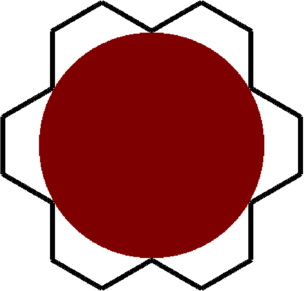}~
\includegraphics[width=0.2\textwidth]{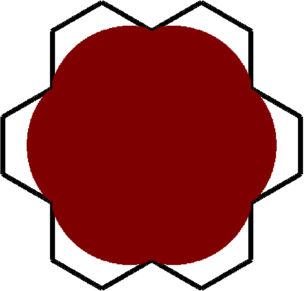}~
\includegraphics[width=0.2\textwidth]{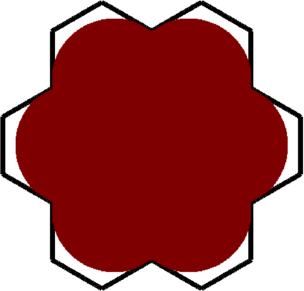}~
\includegraphics[width=0.2\textwidth]{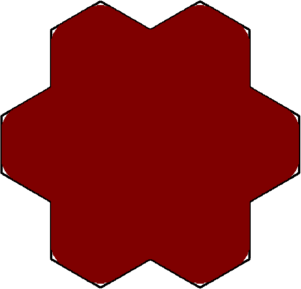}
\caption{The $\alpha$-Cheeger set for a non-convex set in $2D$, for  $\alpha \in \{0.5001,0.75,1,2\}$.}
\label{acheeger_2D}
\end{figure}

\begin{figure}[ht]
\includegraphics[width=0.2\textwidth]{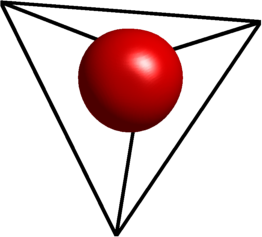}~
\includegraphics[width=0.2\textwidth]{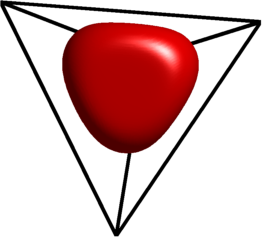}~
\includegraphics[width=0.2\textwidth]{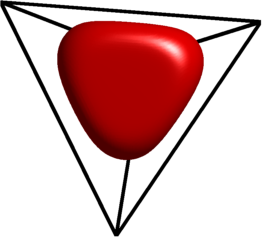}~
\includegraphics[width=0.2\textwidth]{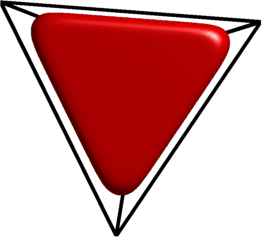}
\caption{The $\alpha$-Cheeger set for a regular tetrahedron in $3D$, for $\alpha \in \{ 0.667,0.9,1,2\}$.}
\label{acheeger_tetra}
\end{figure}

$\bullet$ {\bf Computation of $\alpha$-Cheeger sets.}
In this case, corresponding to $k=1$, there is no need to use the penalization term. We optimize directly the non-penalized ratio between the Modica-Mortola ratio and the volume term with constraints $0\leq u \leq 1$. Examples can be seen in Figure \ref{acheeger_2D} for a domain in $\Bbb{R}^2$ and in Figure \ref{acheeger_tetra} for a domain in $\Bbb{R}^3$. 

\begin{figure}
\begin{tabular}{cccc}
\includegraphics[width=0.18\textwidth]{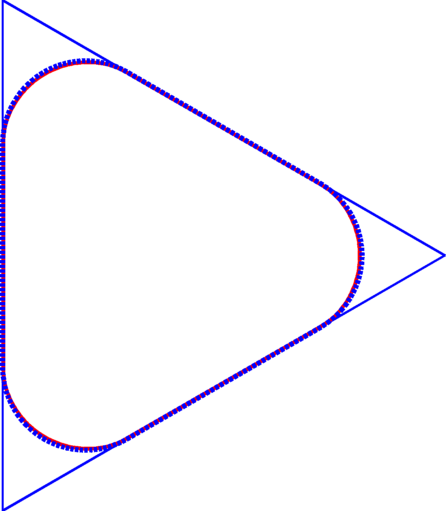} &
\includegraphics[width=0.20\textwidth]{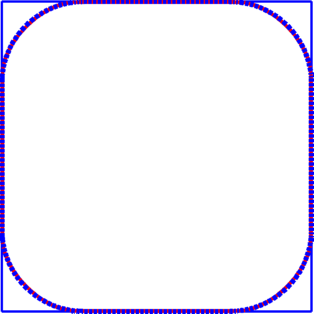} &
\includegraphics[width=0.24\textwidth]{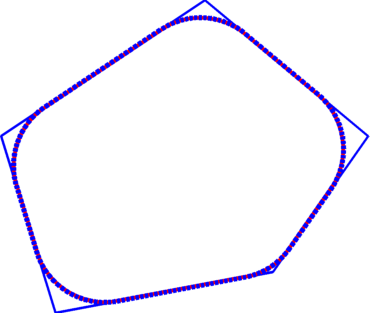} &
\includegraphics[width=0.24\textwidth]{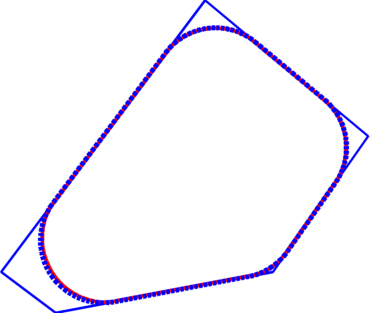}\\
$0.0040$ & $0.0030$ & $0.0030$ & $0.0037$
\end{tabular}
\caption{Comparison between results obtained when minimizing \eqref{optim_func} (red) and the  Kawohl $\&$ Lachand-Robert formula (dotted-blue). Relative errors for the Cheeger constants are also given.}
\label{compar_KLR}
\end{figure}

In order to test the accuracy of our method we compare our algorithm with an implementation of the  Kawohl $\&$ Lachand-Robert explicit formula for finding Cheeger sets $(\alpha=1)$ associated to convex sets in 2D \cite{KLR06}. As can be seen in Figure \ref{compar_KLR} the relaxation algorithm we propose is quite precise. We represent with red the $\varepsilon$-level set of the result obtained when minimizing \ref{optim_func} and with dotted blue the result obtained using the algorithm described in \cite{KLR06}. We justify the choice of the $\varepsilon$-level set by the fact that we impose the boundary condition $u=0$ on the density. Choosing a level set corresponding to a larger value would correspond to a contour which does not touch the boundary of $D$, contrary to the known behaviour of Cheeger sets. In the test cases presented below the results given by the two algorithms are almost indistinguishable. The relative errors, i.e. $|h(\omega_{ap})-h(\omega)|/h(\omega)$ are also provided. Here $h(\omega)$ is the Cheeger constant of $\omega$ and $\omega_{ap}$ is the approximated Cheeger set obtained when minimizing \ref{optim_func}. The computed errors are small in view of the fact that the Cheeger sets are computed with a relaxed formulation. As expected, working on finer meshes leads to better approximations both of the Cheeger sets and of the Cheeger constants.

\begin{figure}[ht]
\includegraphics[width = 0.2 \textwidth]{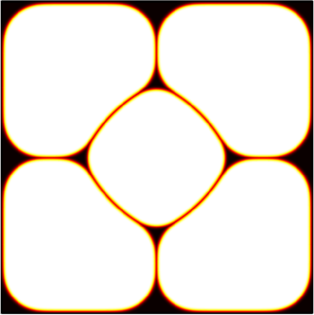}\quad
\includegraphics[width = 0.2 \textwidth]{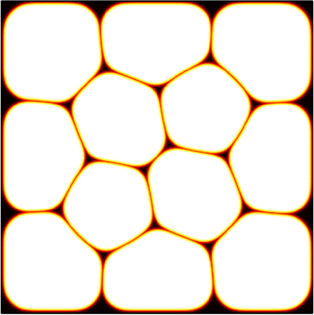}\quad
\includegraphics[width = 0.2 \textwidth]{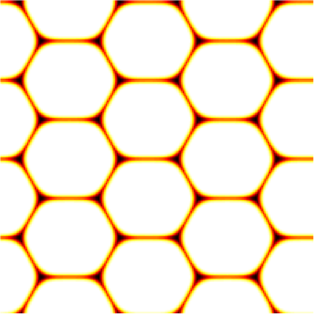}
\caption{Cheeger clusters  for problem \eqref{f:pb2.b} in a square: 5 cells, 12 cells, 16 cells (periodical).}
\label{cheeger_cluster}
\end{figure}

 $\bullet$ {\bf Computation of $\alpha$-Cheeger clusters.} 
Some examples of Cheeger clusters can be seen in Figure \ref{cheeger_cluster}. One can notice immediately that the cells are not necessarily convex, for instance when  $D$ is a square and $n = 5$. The results in the periodic case are in accordance with results in \cite{bfvv17,bf17R}.

 $\bullet$ {\bf Computation of optimal packings.} 
To this aim, we exploit Theorem \ref{t:approx1} and we compute $\alpha$-Cheeger clusters for $\alpha$ very close to $\frac{N-1}{N}$ and $p$ very ``large", which in our computations means at most $100$. Choosing the parameter $\alpha$ close to $\frac{N-1}{N}$ forces the cells in the optimal configurations to be close to disks. We choose to use a $p$-norm approach since this regularizes the non-smooth problem of minimizing the maximal radius of a family of disks. The minimization of a $p$-norm instead of the $\infty$-norm is a natural idea, already used in \cite{BoBN16} for the study of partitions of a domain which minimize the largest fundamental eigenvalue of the Dirichlet-Laplace operator. 

We want to be able to quantify our results, so we use a local refinement procedure as a post-treatment. At the end of the optimization process we have access to the density functions associated to the $\alpha$-Cheeger cells, which are approximately smoothed characteristic functions of disks. Using these density functions we approximate the precise location of the centers of the disks by computing the barycenter of the $0.5$ level-set of each density function. Then we use a very basic local-optimization routine in order to get a precise description of the circle packing that can be compared with existing results in the literature. The algorithm just computes the pairwise distances between centers of the disks and between the centers and the boundary and uses Matlab's algorithm \texttt{fmincon} with the option \texttt{active-set} to perform a gradient-free local optimization of the current configuration. We note that the refinement algorithm is not at all adapted for solving {\it alone} the problem, given random starting points for the centers. It is only useful for locally optimizing the circle packing configuration once localization is obtained by our approximation procedure.

\begin{figure}[ht]
\includegraphics[width=0.2\textwidth]{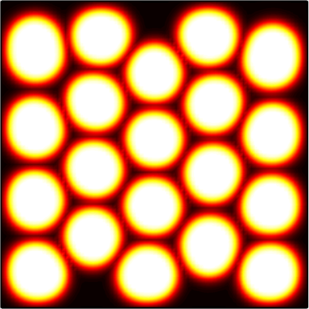}~
\includegraphics[width=0.2\textwidth]{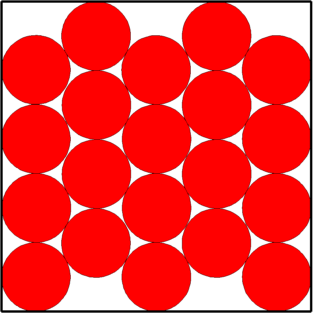}\quad
\includegraphics[width=0.2\textwidth]{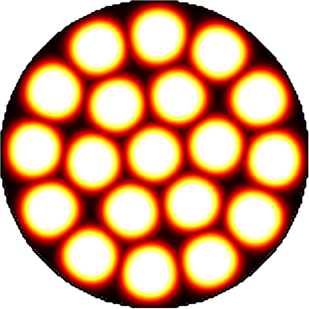}~
\includegraphics[width=0.2\textwidth]{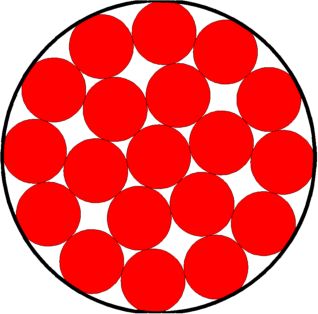}
\vspace{0.1cm}

\bigskip
\includegraphics[width=0.2\textwidth]{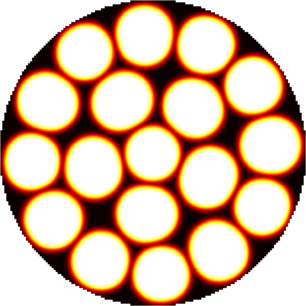}~
\includegraphics[width=0.2\textwidth]{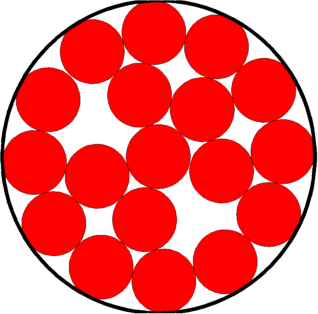}\quad
\includegraphics[width=0.2\textwidth]{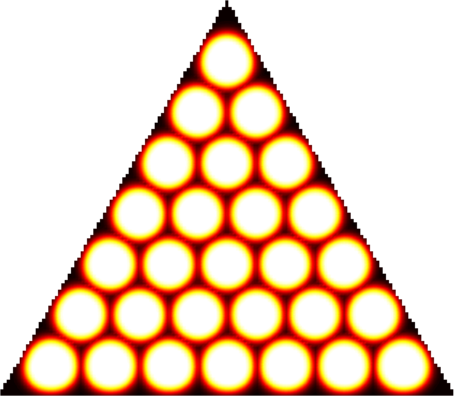}~
\includegraphics[width=0.2\textwidth]{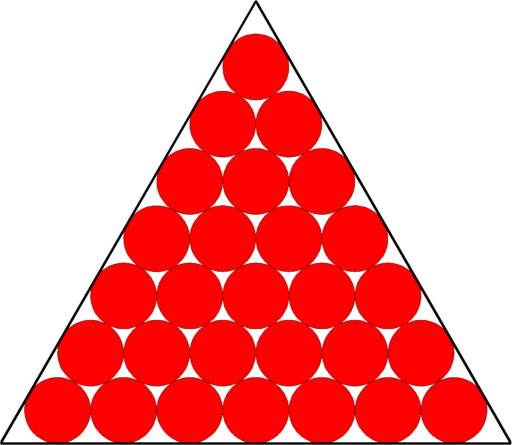}
\caption{Circle packing examples in 2D for problem \eqref{op1}: density representation and local optimization.}
\label{cpack2D}
\end{figure}

It is possible to apply this algorithm for computing both optimal circle packings for domains in $\R^2$ or optimal sphere packings for domains in $\Bbb{R}^3$. In the planar case, we present some computational results in Figure \ref{cpack2D}. In our test cases the numerical algorithm based on the $\Gamma$-convergence result combined with the post-treatment algorithm generally produce configurations which are comparable to the best known results in the literature. We recall that one of the first papers regarding the circle packings in a circle was authored by Kravitz in 1967 \cite{kravitz67}. In this paper we can find a conjecture regarding the $19$-circle packing in a disk. The optimality of this $19$-packing, presented in Figure \ref{cpack2D}, was proved by Fodor in \cite{fodor19}. Extensive numerical results up to thousands of circles were performed and collected on the website \href{http://www.packomania.com/}{\nolinkurl{http://www.packomania.com/}}, maintained by Eckard Specht. In all cases, we compared our results with best ones available, listed on the above cited website. The numerical algorithm manages to capture the right results in cases where the optimal circle packing configuration is unique and rigid, like the case of $19$ disks in a circle or $28$ disks in an equilateral triangle. Moreover, we are able to capture the best known results even in cases where the solution is not unique, like in the case when we have $18$ disks in a circle. 

\bigskip

\bigskip

\begin{figure}[ht]
\includegraphics[width=0.2\textwidth]{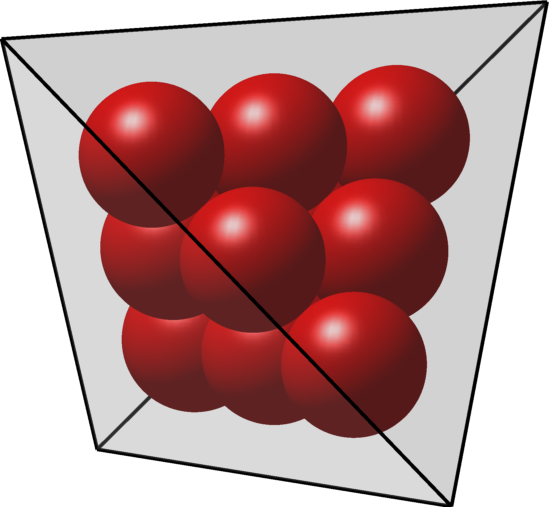}~
\includegraphics[width=0.2\textwidth]{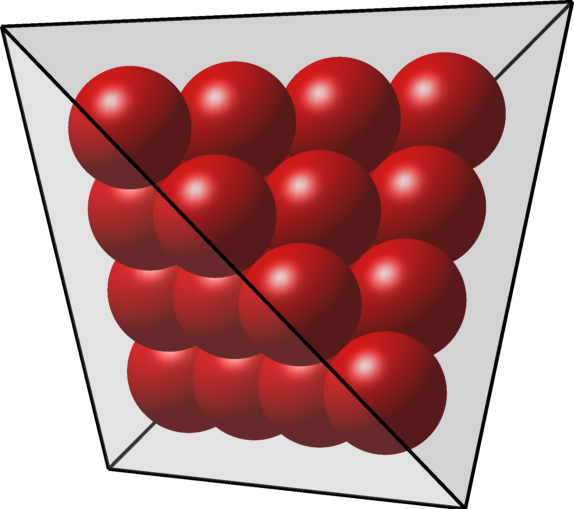}~
\includegraphics[width=0.2\textwidth]{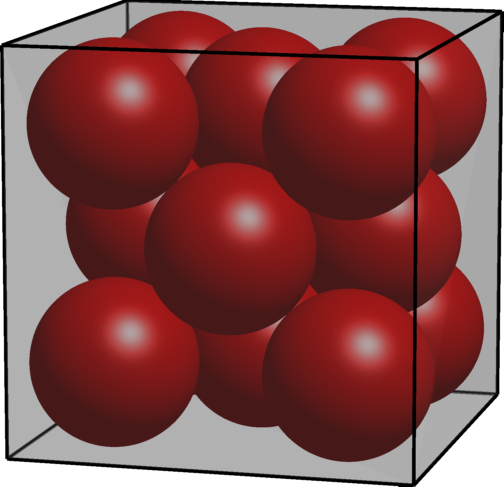}~
\includegraphics[width=0.2\textwidth]{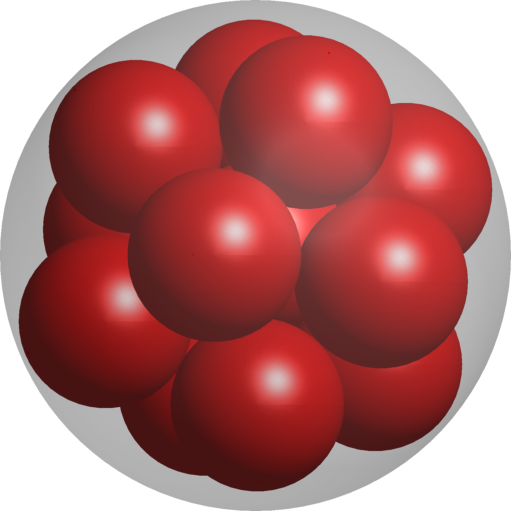}
\caption{Sphere packing examples  in 3D for problem \eqref{op1}.}
\label{spack3D}
\end{figure}

Some examples of computations of optimal spherical packings for domains in $\Bbb{R}^3$ are presented in Figure \ref{spack3D}. In this case, we observe  again a good convergence to the best known configurations.

\bibliographystyle{mybst}
\bibliography{References}

\def\cprime{$'$}
\providecommand{\bysame}{\leavevmode\hbox to3em{\hrulefill}\thinspace}
\providecommand{\MR}{\relax\ifhmode\unskip\space\fi MR }
\providecommand{\MRhref}[2]{%
  \href{http://www.ams.org/mathscinet-getitem?mr=#1}{#2}
}
\providecommand{\href}[2]{#2}
\begin{thebibliography}{10}

\bibitem{BA90}
{S.} Baldo, \emph{Minimal interface criterion for phase transitions in mixtures
  of {C}ahn-{H}illiard fluids}, Ann. Inst. H. Poincar\'e Anal. Non Lin\'eaire
  \textbf{7} (1990), no.~2, 67--90. \MR{1051228}

\bibitem{BogLS17}
{B.} Bogosel, \emph{Efficient algorithm for optimizing spectral partitions},
  Preprint arXiv:1705.08739 (2017).

\bibitem{BoBN16}
{B.} Bogosel and {V.} Bonnaillie-N\"{o}el, \emph{Minimal Partitions for p-norms
  of Eigenvalues}, Preprint arXiv:1612.07296 (2016).

\bibitem{BO16}
{B.} Bogosel and {\'E.}~Oudet, \emph{Qualitative and numerical analysis of a
  spectral problem with perimeter constraint}, SIAM J. Control Optim.
  \textbf{54} (2016), no.~1, 317--340. \MR{3459973}

\bibitem{BBO09}
{B.} Bourdin, {D.} Bucur, and {\'E.}~Oudet, \emph{Optimal partitions for
  eigenvalues}, SIAM J. Sci. Comput. \textbf{31} (2009/10), no.~6, 4100--4114.
  \MR{2566585}

\bibitem{bf17R}
{D.} Bucur and {I.} Fragal{\`a}, \emph{On the honeycomb conjecture for {R}obin
  {L}aplacian eigenvalues}, preprint CVGMT (2017).

\bibitem{bfvv17}
{D.} Bucur, {I.} Fragal\`a, {B.} Velichkov, and {G.} Verzini, \emph{On the
  honeycomb conjecture for a class of minimal convex partitions}, Arxiv
  Preprint, arXiv:1703.05383 (2017).

\bibitem{CaffLin}
{L. A.} Caffarelli and {F. H.} Lin, \emph{An optimal partition problem for
  eigenvalues}, J. Sci. Comput. \textbf{31} (2007), no.~1-2, 5--18.

\bibitem{CCP09}
{G.} Carlier, {M.} Comte, and {G.} Peyr\'e, \emph{Approximation of maximal
  {C}heeger sets by projection}, M2AN Math. Model. Numer. Anal. \textbf{43}
  (2009), no.~1, 139--150. \MR{2494797}

\bibitem{Car17}
M.~Caroccia, \emph{Cheeger {N}-clusters}, Calc. Var. Partial Differential
  Equations \textbf{56} (2017), no.~2, 56:30. \MR{3610172}

\bibitem{cfm09}
{V.} Caselles, {G.} Facciolo, and {E.} Meinhardt, \emph{Anisotropic {C}heeger
  sets and applications}, SIAM J. Imaging Sci. \textbf{2} (2009), no.~4,
  1211--1254. \MR{2559165}

\bibitem{DM}
{G.} Dal~Maso, An introduction to {$\Gamma$}-convergence, Progress in Nonlinear
  Differential Equations and their Applications, vol.~8, Birkh\"auser Boston,
  Inc., Boston, MA, 1993.

\bibitem{fodor19}
{F.} Fodor, \emph{The densest packing of {$19$} congruent circles in a circle},
  Geom. Dedicata \textbf{74} (1999), no.~2, 139--145. \MR{1674049}

\bibitem{cpack10}
A.~Grosso, A.~R. M. J.~U. Jamali, M.~Locatelli, and F.~Schoen, \emph{Solving
  the problem of packing equal and unequal circles in a circular container}, J.
  Global Optim. \textbf{47} (2010), no.~1, 63--81. \MR{2609042}

\bibitem{Hales}
{T. C.} Hales, \emph{The honeycomb conjecture}, Discrete Comput. Geom.
  \textbf{25} (2001), no.~1, 1--22.

\bibitem{KLR06}
{B.} Kawohl and {T.} Lachand-Robert, \emph{Characterization of {C}heeger sets
  for convex subsets of the plane}, Pacific J. Math. \textbf{225} (2006),
  no.~1, 103--118. \MR{2233727}

\bibitem{kravitz67}
{S.} Kravitz, \emph{Packing {C}ylinders into {C}ylindrical {C}ontainers}, Math.
  Mag. \textbf{40} (1967), no.~2, 65--71. \MR{1571666}

\bibitem{LRO05}
{T.} Lachand-Robert and {E.} Oudet, \emph{Minimizing within convex bodies using
  a convex hull method}, SIAM J. Optim. \textbf{16} (2005), no.~2, 368--379.
  \MR{2197985}

\bibitem{Leo}
{G. P.} Leonardi, \emph{An overview on the Cheeger problem}, Pratelli, {A.},
  Leugering, {G.} (eds.) New trends in shape optimization., International
  Series of Numerical Mathematics, Springer (Switzerland), vol. 166, 2016,
  pp.~117--139.

\bibitem{cpackLB11}
C.~O. L\'opez and J.~E. Beasley, \emph{A heuristic for the circle packing
  problem with a variety of containers}, European J. Oper. Res. \textbf{214}
  (2011), no.~3, 512--525. \MR{2820172}

\bibitem{modica-mortola}
{L.} Modica and {S.} Mortola, \emph{Un esempio di {$\Gamma ^{-}$}-convergenza},
  Boll. Un. Mat. Ital. B (5) \textbf{14} (1977), no.~1, 285--299. \MR{0445362
  (56 \#3704)}

\bibitem{Oudet11}
{E.} Oudet, \emph{Approximation of partitions of least perimeter by
  {$\Gamma$}-convergence: around {K}elvin's conjecture}, Exp. Math. \textbf{20}
  (2011), no.~3, 260--270. \MR{2836251}

\bibitem{Pa}
{E.} Parini, \emph{An introduction to the {C}heeger problem}, Surv. Math. Appl.
  \textbf{6} (2011), 9--21.

\bibitem{PS17}
{A.} Pratelli and {G.} Saracco, \emph{On the generalized {C}heeger problem and
  an application to 2d strips}, Rev. Mat. Iberoam. \textbf{33} (2017), no.~1,
  219--237. \MR{3615449}

\bibitem{lbfgs}
{L.} Stewart, \emph{MATLAB LBFGS Wrapper},
  \url{http://www.cs.toronto.edu/~liam/software.shtml}.

\end{thebibliography}

\end{document}